\documentclass[12pt,leqno]{article}
\usepackage{amsmath, amsthm, amsfonts, amssymb, color}
\usepackage{mathrsfs}
\theoremstyle{definition}

\newcommand{\scr}[1]{\mathscr #1}
\definecolor{wco}{rgb}{0.5,0.2,0.3}

\numberwithin{equation}{section} \theoremstyle{remark}

\newcommand{\ua}{\uparrow}

\title{{\bf A Sharp Liouville
Theorem for Elliptic Operators}
 \footnote{Supported in
 part by NNSFC(10721091) and the 973-Project
  and by the Italian  M.I.U.R. research project
``Equazioni di Kolmogorov''.} }

\author{ Enrico Priola \footnote{
 enrico.priola@unito.it}
 \\
{\footnotesize  {Department of Mathematics, University of Torino,
 Torino, 10123, Italy} }
\\ Feng-Yu Wang \footnote{
 wangfy@bnu.edu.cn;
F.Y.Wang@swansea.ac.uk }
\\ \footnotesize
{School of Mathematical Sci. and Lab. Math. Com. Sys.,
Beijing Normal
University, Beijing 100875, China}\\
\footnotesize{and}\\ \footnotesize{Department of Mathematics,
Swansea University, Singleton Park, SA2 8PP, UK} }

\begin{document}

\def\R{\mathbb R}  \def\ff{\frac} \def\ss{\sqrt} \def\BB{\mathbb
B}
\def\N{\mathbb N} \def\kk{\kappa} \def\m{{\bf m}}
\def\dd{\delta} \def\DD{\Delta} \def\vv{\varepsilon} \def\rr{\rho}
\def\<{\langle} \def\>{\rangle} \def\GG{\Gamma} \def\gg{\gamma}
  \def\nn{\nabla} \def\pp{\partial} \def\tt{\tilde}
\def\d{\text{\rm{d}}} \def\bb{\beta} \def\aa{\alpha} \def\D{\scr D}
\def\EE{\mathbb E} \def\si{\sigma} \def\ess{\text{\rm{ess}}}
\def\beg{\begin} \def\beq{\begin{equation}}  \def\F{\scr F}
\def\Ric{\text{\rm{Ric}}} \def\Hess{\text{\rm{Hess}}}\def\B{\scr B}
\def\e{\text{\rm{e}}} \def\ua{\underline a} \def\OO{\Omega}  \def\oo{\omega}     \def\tt{\tilde} \def\Ric{\text{\rm{Ric}}}
\def\cut{\text{\rm{cut}}} \def\P{\mathbb P} \def\ifn{I_n(f^{\bigotimes n})}
\def\C{\scr C}      \def\aaa{\mathbf{r}}     \def\r{r}
\def\gap{\text{\rm{gap}}} \def\prr{\pi_{{\bf m},\varrho}}  \def\r{\mathbf r}
\def\Z{\mathbb Z} \def\vrr{\varrho} \def\ll{\lambda}
\def\L{\scr L}\def\Tt{\tt} \def\TT{\tt}\def\II{\mathbb I}
\def\i{{\rm i}}\def\Sect{{\rm Sect}}\def\E{\mathbb E}\def\LL{\Lambda}

\maketitle
\begin{abstract} We introduce a new condition
 on
 elliptic operators $L= \frac{1}{2}\triangle + b \cdot \nabla$
  which ensures the validity of the
  Liouville property for
 bounded solutions  to  $Lu=0$ on $\R^d$. Such condition is
 sharp when $d=1$. We extend  our Liouville theorem  to more
 general second order operators in non-divergence form
  assuming a Cordes type condition.
 \end{abstract}


\section{Introduction}
Let
 $$ L = \frac{1}{2}\sum_{i,j =1}^d q_{ij}(x) D_{ij} +
 \sum_{i=1}^d b_i (x) D_i $$
 be a uniformly elliptic second order differential operator
  on $\R^d$ with continuous coefficients $q_{ij}$ and $b_i$
 (here $D_{ij} = \ff{\pp^2}{\pp x_i\pp x_j}$
 and $D_i= \ff{\pp}{\pp x_i}, \ 1\le i,j\le d$).
 Recall that
a smooth real function $u$ on $\R^d$ is called $L$-harmonic if $Lu=0$ holds on
$\R^d.$ An operator $L$ is said to possess the Liouville property
 when all bounded $L$-harmonic functions are constant (or,
 equivalently,
  when a two-sided Liouville theorem holds for $L$).
  Such property is also of interest
   for the study of non-linear
 PDEs of the form $\triangle u + F(u) =0$ (see
    e.g. \cite{Ba, BCN}).


There are a plenty of results on the Liouville property.
 Let $\lambda_0>0$ be the ellipticity constant of $L$. A typical
condition implying the Liouville property  is the following (see
e.g. \cite{BF, PW, PZ}):
\begin{align} \label{qx}
 \frac{1}{2 \lambda_0} \| q(x) - q(x +h) \|^2 + 2
 \<  b(x+h) - b (x), h \> \le 0, \;\; x, \, h \in \R^d
\end{align}
 (given a $d \times d$
 real matrix $A $, we denote by $\| A\|$
its Hilbert-Schmidt norm; moreover
 $\< \cdot, \cdot \>$ is the Euclidean inner product
  in $\R^d$). However this is not   completely  satisfactory for two
  reasons. The first one is that when
   $b(x)$ is constant the matrix
   $q(x)$ must be constant as well. This is a restriction
    since      it is known  that
   the Liouville property holds when $b$ is constant
   and $q(x) $ is variable (see \cite[Corollary 4.1]{KS}
   which is  a consequence of  invariant Harnack
   inequalities).

 The second weak point of \eqref{qx} is that  when
 $q(x)$ is the identity, i.e., we are considering
  $L_0=\ff 1 2 \DD+b\cdot\nn$,
  such hypothesis is not optimal even
  when
  $d=1$.
 The aim of this note is to find out a sharp and easy to check
 criterion ensuring  the Liouville  property   for
 $L_0$. Our condition is sharp when  $d=1$; indeed   if
 this does not hold one can construct counterexamples
  of operators $L_0$
  without the Liouville property.

  We prove our  Liouville type
   theorem
    in the  more general setting of
  elliptic operators  $L$, with $q(x)$ variable,
   imposing an additional Cordes type condition (see \cite{Cor}).


\vskip 2mm
 To explain the motivation of our
desired condition for the Liouville property,
 let us start with a
  one-dimensional
 example
$$
L_0=\ff 1 2 \ff{\d^2}{\d x^2} + \ff{x}{2+x^2}\Big(\dd+\ff{2}{\log(2+
x^2)}\Big)\ff{\d}{\d x},$$
where $\dd$ is a constant.  It is easy to
see that a harmonic function of $L_0$ has the form
$$u(x)= c_1 +c_2 \int_0^x \ff {\d r}
{(2+r^2)^{\dd}\log^2(2+r^2)},\ \ \  x\in \R,$$ where $c_1, c_2$ are
constants. Thus, all bounded harmonic functions are constant if any
only if $\dd< 1/2.$  In order to reduce this condition to a usual
monotonicity condition on the drift $b(x)=
\ff{x}{2+x^2}\big(\dd+\ff{2}{\log(2+ x^2)}\big),$ we note that
 (using also that $b$ is odd)
 $$\lim_{s\to\infty}\sup_{|x-y|=s}
 (x-y) \big( b(x)-b(y) \big)=
  \lim_{s\to\infty}  s( b(s/2)-b(-s/2))= 4\dd.$$
  Then the statement can be reformulated as
 all bounded $L_0$-harmonic functions are constant if and only if
 $$
 \lim_{s\to\infty}\sup_{|x-y|=s} (x-y)  \big( b(x)-b(y) \big) \,<\,
 2.$$
 In general, let e.g. $L_0=\ff 1 2 \DD+b\cdot\nn$ on $\R^d$, we may wish to prove the Liouville
 property of $L_0$ under the following hypothesis
\beq\label{1.12} \limsup_{s\to\infty}\sup_{|x-y|=s} \<x-y,
b(x)-b(y)\> \, < \, 2.
\end{equation}
 This follows immediately from our main result.

\section {Main theorem}

 We prove  a Liouville type theorem for bounded
 space-time harmonic functions.
 Recall that a smooth function $u$
  on $[0,\infty)\times \R^d$ is called
 space-time harmonic for $L$, if $\pp_t u+Lu=0$ holds.
 To state our main result, we make the following assumptions.

\paragraph{(H)}
$(i)$ \ The coefficients $b(x)$ and  $q(x) $
 are  continuous, and, for any  $\ll>0$,
 $\omega (s): = \sup_{|x-y|\le s} \{\ll \| q(x) - q(y) \|^2+2\<x-y, b(x)-b(y)\>\}$ satisfies
 $$
 \int_0^1 \frac{\omega(s)}{s}ds < \infty;
 $$
$(ii)$  there exist two constants $0 < \lambda_0 < \Lambda_0$ such that
  $$ \lambda_0 |h|^2 \le \sum_{i,j =1}^n q_{ij}(x) h_i h_j \le
\Lambda_0 |h|^2,\;\; x, \; h \in \R^d.$$

\beg{thm} \label{T1.1}  Assume {\bf (H)}. If

\beq\label{1.1} \limsup_{s\to\infty}\sup_{|x-y|=s} \<x-y,
b(x)-b(y)\> <2\ll_0 - \ff d 2(\LL_0-\ll_0),\end{equation} then any
bounded
 space-time harmonic function for $L$ is constant.  \end{thm}

\beg{proof}
  We will suitably apply \cite[Theorem 3.6]{PW}.
 To this purpose,
   we have to consider  a  coupling for $L$.
 By (\ref{1.1}) we may take constants $\mu, s_0>0$ and $s_1\in \R$ such that $\mu < \lambda_0$ and

 \beq\label{1.2}\sup_{|x-y|=s} \<x-y, b(x)-b(y)\> \le s_1<2\mu -\ff 1 2  d(\LL_0-\mu),\ \ \ s\ge s_0.\end{equation}
 Define a symmetric positive definite matrix $\sigma (x)$, such that $\sigma (x)^2 + \mu I = q(x)$, $x \in \R^d$.
Clearly we have $\sigma ^2(x) \ge (\lambda_0 - \mu )I $. We construct a coupling as in Section 3.1 of \cite{PW},
replacing the ellipticity constant $\lambda_0$ with $\mu$
 (note that under
  our assumptions the associated diffusion process does not
  explode).
 Applying   \cite[Lemma 3.3]{PW} we deduce that
$$
\|  \sigma (x) - \sigma (y) \|^2 \le \frac{1}{4 (\lambda_0 - \mu )}
\| q (x) - q (y) \|^2,\;\;\; x,\, y \in \R^d.$$ Combining this  with
{\bf (H)}(i) for $\ll= \ff 1 {4(\ll_0-\mu)}, $
 we obtain

\begin{align}\label{1.3}
& \|\si(x)-\si(y)\|^2+ 2\<x-y, b(x)-b(y)\> \le \oo(|x-y|)\ \;\;\;
\text{for}\ x,y\in \R^d,\  \text{and}
\\ \nonumber & \int_0^{s_0}
\ff{\oo(s)} s\d s<\infty.
\end{align}
On the other hand, since $\sigma (x)^2 \le (\Lambda_0 - \mu)I$, we
have  $0 \le
 \sigma (x) \le (\Lambda_0 - \mu)^{1/2}I$, for any $x \in \R^d$.
  Thus
 $$
- (\Lambda_0 - \mu)^{1/2}I \le \sigma (x) - \sigma (y) \le
(\Lambda_0 - \mu)^{1/2}I, \;\; x, y \in \R^d.
$$
 We deduce that $ 0 \le \big( \sigma (x) - \sigma (y) \big)^2 \le (\Lambda_0 -
 \mu)I$ and so
 $$ \|  \sigma (x) - \sigma (y) \|^2
 = \text{Tr} \big[( \sigma (x) - \sigma (y)  )^2 \big] \le
 d (\Lambda_0 - \mu),\;\; x, y \in \R^d. $$ Combining this with (\ref{1.2}) we obtain

 $$\|\si(x)-\si(y)\|^2+2\<x-y, b(x)-b(y)\>\le 2s_1+ d(\LL_0-\mu)=:s_2<4\mu,\ \ \ |x-y|\ge s_0.$$
From this and (\ref{1.3}) we conclude that

$$\|\si(x)-\si(y)\|^2+2\<x-y, b(x)-b(y)\>\le |x-y|g(|x-y|),\ \ \ x,y\in \R^d$$ holds for
 $$
 g (s) : = \ff{\oo(s)}
 s 1_{[0,s_0]}(s) + \frac{s_2}{ s}1_{(s_0,\infty)},\ \ \  s >0.
 $$ Since by {\bf (H)}

 $$c:=\int_0^{s_0} g(s)\d s<\infty,$$ we have

 \beg{equation*} \beg{split}
 & \int_0^{\infty} \exp \Big( -\frac{1}{4 \mu} \int_0^r g(s)ds \Big ) \, dr \\
 &\ge
 \int_1^{\infty}   \exp \Big( -\frac{1}{4 \mu} \int_0^{s_0} g(s)ds \Big ) \,
  \exp \Big( -\frac{1}{4 \mu} \int_{s_0}^r g(s)ds \Big ) \, dr
\\
&\ge  \e^{-c_1}  \int_{1}^{\infty}  s^{-s_2/[4\mu]}\d s=
  \infty\end{split}\end{equation*}
 since $s_2< 4\mu$. Applying
   \cite[Theorem 3.6]{PW}, we get the assertion.
\end{proof}


\begin{thebibliography}{pippo1}
 \footnotesize{


\bibitem{Ba} Barlow, M. T., On the Liouville property for divergence
 form operators.  Canad. J. Math. 50 (1998), no. 3,
 487-496.

\bibitem{BCN} Berestycki, H., Caffarelli, L., Nirenberg, L., Further
qualitative properties for elliptic equations in unbounded domains.
Dedicated to Ennio De Giorgi.  Ann. Scuola Norm. Sup. Pisa Cl. Sci. (4)  25  (1997),  no. 1-2,  (1998) 69-94.


\bibitem{BF} Bertoldi M., Fornaro S,
 Gradient
 estimates in parabolic problems with unbounded coefficients. Studia
 Math. 165 (2004), no. 3, 221-254.



\bibitem{Cor} Cordes, H. O.
 { Uber die
 erste Randwertaufgabe bei quasilinearen Differentialgleichungen
zweiter Ordnung in mehr als zwei Variablen.} (German) Math. Ann. 131
(1956), 278-312.




\bibitem{KS}
 Krylov, N. V.,  Safonov, M. V.,  A property of the solutions of
 parabolic equations with measurable coefficients. Izv. Akad. Nauk
 SSSR Ser. Mat. 44 (1980), no. 1, 161-175.

\bibitem{PW} Priola, E., Wang, Feng-Yu, {
   Gradient estimates for diffusion
 semigroups with singular coefficients.}
  J. Funct. Anal.  236  (2006),  244-264.



\bibitem{PZ}  Priola, E., Zabczyk, J.,
   Liouville theorems for nonlocal
operators. { J. Funct. Anal.,} 216 (2004),
 455-490.





}
\end{thebibliography}
\end{document}